\documentclass[12pt]{article}
\usepackage{amsfonts}
\usepackage{latexsym,amsmath,amssymb}
\usepackage{graphics,epstopdf}
\usepackage[pdftex]{graphicx}
\usepackage{color}
\usepackage{hyperref}
\hypersetup{colorlinks=false}
\usepackage[ top=1.5cm, bottom=1.5cm, outer=2cm, inner=2cm,
heightrounded, marginparwidth=2.5cm, marginparsep=2cm] {geometry}
\numberwithin{equation}{section}
\begin{document}

\bigskip

\begin{center}
{\large \textbf{Stability of Some Positive Linear Operators on Compact Disk}}%
\\[0pt]
\vspace{.7cm} \textbf{M. Mursaleen}, \textbf{Khursheed J. Ansari} and
\textbf{Asif Khan} \bigskip

\vspace{.2cm} Department of Mathematics, Aligarh Muslim University, Aligarh
202002, India\\[0pt]
Email: mursaleenm@gmail.com; ansari.jkhursheed@gmail.com; asifjnu07@gmail.com%
\\[0pt]
\end{center}

\bigskip

{\footnotesize {\textbf{Abstract.}} Recently, Popa and Ra\c{s}a \cite{popa2,
popa3} have shown the (in)stability of some classical operators defined on $%
[0,1]$ and found best constant when the positive linear operators are stable
in the sense of Hyers-Ulam. In this paper we show Hyers-Ulam (in)stability
of complex Bernstein-Schurer operators, complex Kantrovich-Schurer operators
and Lorentz operators on compact disk. In the case when the operator is
stable in the sense of Hyers and Ulam, we find the infimum of Hyers-Ulam
stability constants for respective operators.}

\bigskip

\bigskip

\noindent {\large \textbf{1. Introduction}}\newline
\noindent The equation of homomorphism is stable if every \textquotedblleft
approximate\textquotedblright\ solution can be approximated by a solution of
this equation. The problem of stability of a functional equation was
formulated by S.M. Ulam \cite{ulam} in a conference at Wisconsin University,
Madison in 1940: \textquotedblleft Given a metric group $(G,.,\rho )$, a
number $\varepsilon >0$ and a mapping $f:G\rightarrow G$ which satisfies the
inequality $\rho (f(xy),f(x)f(y))<\varepsilon $ for all $x,y\in G$, does
there exist a homomorphism $a$ of $G$ and a constant $k>0$, depending only
on $G$, such that $\rho (a(x),f(x))\leq k\varepsilon $ for all $x\in G?$%
\textquotedblright\ If the answer is affirmative the equation $%
a(xy)=a(x)a(y) $ of the homomorphism is called stable; see \cite{ras,isac}.
The first answer to Ulam's problem was given by D.H. Hyers \cite{hyer} in
1941 for the Cauchy functional equation in Banach spaces, more precisely he
proved: \textquotedblleft Let $X,~Y$ be Banach spaces, $\varepsilon $ a
non-negative number, $f:X\rightarrow Y$ a function satisfying $\Vert
f(x+y)-f(x)-f(y)\Vert \leq \varepsilon $ for all $x,y\in X$, then there
exists a unique additive function with the property $\Vert f(x)-a(x)\Vert
\leq \varepsilon $ for all $x\in X$.\textquotedblright\ Due to the question
of Ulam and the result of Hyers this type of stability is called today
Hyers-Ulam stability of functional equations. A similar problem was
formulated and solved earlier by G. P\'{o}lya and G. Szeg\"{o} in \cite{pol}
for functions defined on the set of positive integers. After Hyers result a
large amount of literature was devoted to study Hyers-Ulam stability for
various equations. A new type of stability for functional equations was
introduced by T. Aoki \cite{aok} and Th.M. Rassias \cite{ras0} by replacing $%
\varepsilon $ in the Hyers theorem with a function depending on $x$ and $y$,
such that the Cauchy difference can be unbounded. For other results on the
Hyers-Ulam stability of functional equations one can refer to \cite{ras,
ansari,gajda,Brz,bri,rus,urs,sin,rus1,bota,bota1,pet,brz1,ansari1}.
 The Hyers-Ulam stability of linear operators was considered for the
first time in the papers by Miura, Takahasi et al. (see \cite{hat,miu,tak}).
Similar type of results are obtained in \cite{miu1} for weighted composition
operators on $C(X)$, where $X$ is a compact Hausdorff space. A result on the
stability of a linear composition operator of the second order was given by
J. Brzdek and S.M. Jung in \cite{jung}.

\parindent=8mm Recently, Popa and Ra\c{s}a obtained \cite{rasa2} a result on
Hyers-Ulam stability of the Bernstein-Schnabl operators using a new approach
to the Fr\'{e}chet functional equation, and in \cite{popa2, popa3}, they
have shown the (in)stability of some classical operators defined on $[0,1]$
and found the best constant for the positive linear operators in the sense
of Hyers-Ulam.

\parindent=8mm Motivated by their work, in this paper, we show the
(in)stability of some complex positive linear operators on compact disk in
the sense of Hyers-Ulam. We find the infimum of the Hyers-Ulam stability
constants for complex Bernstein-Schurer operators and complex
Kantrovich-Schurer operators on compact disk. Further we show that Lorentz
polynomials are unstable in the sense of Hyers-Ulam on a compact disk.%
\newline

\noindent {\large \textbf{2. The Hyers-Ulam stability property of operators}}%
\newline

\noindent In this section, we recall some basic definitions and results on
Hyers-Ulam~stability~property which form the background of our main results.%
\newline

\noindent {\large \textbf{Definition 2.1.}} Let $A$ and $B$ be normed spaces
and $T$ a mapping from $A$ into $B$. We say that $T$ has the \textit{%
Hyers-Ulam~stability~property} (briefly, $T$ is \textit{HU-stable}) \cite%
{miu1} if there exists a constant $K$ such that:\newline
(i) for any $g\in T(A)$, $\varepsilon >0$ and $f\in A$ with $\Vert Tf-g\Vert
\leq \varepsilon $, there exists an $f_{0}\in A$ such that $Tf_{0}=g$ and $%
\Vert f-f_{0}\Vert \leq K\varepsilon $. The number $K$ is called a \textit{%
HUS~constant~of}$~T$ , and the infimum of all HUS constants of $T$ is
denoted by $K_{T}$. Generally, $K_{T}$ is not a HUS constant of T ; see \cite%
{hat,miu}.\newline

\parindent=8mm Let now $T$ be a bounded linear operator with the kernel
denoted by $N(T)$ and the range denoted by $R(T)$. Consider the one-to-one
operator $\widetilde{T}$ from the quotient space $A/N(T)$ into $B$:\newline
\begin{equation*}
\widetilde{T}(f+N(T))=Tf,~~f\in A,\newline
\end{equation*}%
and the inverse operator $\widetilde{T}^{-1}:R(T)\rightarrow A/N(T)$.\newline

\noindent {\large \textbf{Theorem 2.2.}}(\cite{miu1}). \textit{Let $A$ and $%
B $ be Banach spaces and $T:A\rightarrow B$ be a bounded linear operator.
Then the following statements are equivalent:}\newline
(a) $T$ \textit{is HU-stable;}\newline
(b) $R(T)$ \textit{is closed;}\newline
(c) $\widetilde{T}^{-1}$ \textit{is bounded.}\newline
\textit{Moreover, if one of the conditions} (a), (b), (c) \textit{is
satisfied, then }$K_{T}=\Vert \widetilde{T}^{-1}\Vert $.\newline

\noindent {\large \textbf{Remark 2.3.}} (1) Condition (i)
expresses the Hyers-Ulam stability of the equation\newline
$~~~~Tf=g$, where $g\in R(T)$ is given and $f\in A$ is unknown.\newline
(2) If $T:A\rightarrow B$ is a bounded linear operator, then (i) is
equivalent to:\newline
(ii) for any $f\in A$ with $\Vert Tf\Vert \leq 1$ there exists an $f_{0}\in
N(T)$ such that $\Vert f-f_{0}\Vert \leq K$, (see [13]).\newline
So, in what follows, we shall study the HU-stability of a bounded linear
operator\newline
$T:A\rightarrow B$ by checking the existence of a constant $K$ for which
(ii) is satisfied, or equivalently, by checking the boundedness of $%
\widetilde{T}^{-1}$.\newline

\parindent=8mm The main results used in our approach for obtaining, in some
concrete cases, the explicit value of $K_{T}$ are the formula given above
and a result by Lubinsky and Ziegler \cite{lub} concerning coefficient
bounds in the Lorentz representation of a polynomial. Let $p\in \it{\Pi _{n}}$,
 where $\it{\Pi _{n}}$ is the set of all polynomials of
degree at most $n$ with real coefficients. Then $p$ has a unique Lorentz
representation of the form%
\begin{equation*}
p(x)=\sum\limits_{k=0}^{n}c_{k}x^{k}(1-x)^{n-k},~\eqno(2.1)
\end{equation*}%
\newline
where $c_{k}\in \mathbb{R}$, $k=0,1,\cdots ,n$. Remark that, in fact, it is
a representation in Bernstein-B$\acute{e}$zier basis. Let $T_{n}$ denote the
usual $n$th degree Chebyshev polynomial of the first kind. Then the
following representation holds (see \cite{lub}):\newline
\begin{equation*}
T_{n}(2x-1)=\sum\limits_{k=0}^{n}d_{n,k}x^{k}(1-x)^{n-k}(-1)^{n-k},~\eqno%
(2.2)
\end{equation*}%
where
\begin{equation*}
d_{n,k}:=\sum\limits_{j=0}^{\min \{k,n-k\}}\left(
\begin{array}{c}
n \\
2j%
\end{array}%
\right) \left(
\begin{array}{c}
n-2j \\
k-j%
\end{array}%
\right) 4^{j},~~k=0,1,\cdots ,n.
\end{equation*}%
\newline
It is proved in \cite{popa2} that $d_{n,k}=\left(
\begin{array}{c}
2n \\
2k%
\end{array}%
\right) ,~~k=0,1,\cdots ,n$. Therefore
\begin{equation*}
T_{n}(2x-1)=\sum\limits_{k=0}^{n}{\binom{2n}{2k}}(-1)^{n-k}x^{k}(1-x)^{n-k}.%
\newline
\end{equation*}

\noindent {\large \textbf{Theorem 2.4.}}(\emph{Lubinsky and Ziegler \cite%
{lub}). Let $p$ have the representation (2.1), and let $0\leq k\leq n$. Then}
\begin{equation*}
|c_{k}|\leq d_{n,k}.\Vert p\Vert _{\infty }\newline
\end{equation*}%
\emph{with equality if and only if $p$ is a constant multiple of $%
T_{n}(2x-1) $} where $\Vert p\Vert _{\infty }=\max\limits_{x\in[a,b]}|P(x)|$.\newline

\parindent=8mm Let $C[0,1]$ be the space of all continuous, real-valued
functions defined on $[0,1]$, and $C_{b}[0,+\infty )$ the space of all
continuous, bounded, real-valued functions on $[0,+\infty )$. Endowed with
the supremum norm, they are Banach spaces.\newline

\parindent=8mm Popa and Ra\c{s}a have shown the Hyers-Ulam stability of the
following operators:\newline

\noindent {\large \textbf{(i) Bernstein operators}} (\cite{popa2}) \newline

For each integer $n\geq 1$, the sequence of classical Bernstein operators $%
B_{n}:C[0,1]\rightarrow C[0,1]$ is defined by (see [1])%
\begin{equation*}
B_{n}f(x)=\sum\limits_{k=0}^{n}{\binom{n}{k}}x^{k}(1-x)^{n-k}f\biggl{(}\frac{k}{n}%
\biggl{)},~~~f\in C[0,1],~n\geq 1.\newline
\end{equation*}%
It is stable in the Hyers-Ulam sense and the best Hyers-Ulam stability
constant is given by%
\begin{equation*}
~~K_{B_{n}}={\binom{{2n}}{{2[\frac{n}{2}]}}},n\in \mathbb{N}.
\end{equation*}%
\newline

\noindent {\large \textbf{{(ii) Sz\'{a}sz-Mirakjan operators}}} (\cite{popa2}%
)\newline

\noindent The $n$th Sz\'{a}sz-Mirakjan operator $L_{n}:C_{b}[0,+\infty
)\rightarrow C_{b}[0,+\infty )$ defined by (see \cite{cam}, pp. 338)%
\begin{equation*}
~~~~L_{n}f(x)=e^{-nx}\sum\limits_{j=0}^{\infty }f\biggl{(}\frac{j}{n}\biggl{)}\frac{n^{j}}{j!%
}x^{j},~~x\in \lbrack 0,+\infty )
\end{equation*}%
\newline
$~$is not stable in the sense of Hyers and Ulam for each $n\geq 1$.\newline

\noindent {\large \textbf{(iii) Beta operators}} (\cite{popa2})\newline

For each $n\geq 1$, the Beta operator $B_{n}:C[0,1]\rightarrow C[0,1]$
defined by \cite{lup}%
\begin{equation*}
~~~~L_{n}f(x):=\frac{\int_{0}^{1}t^{nx}(1-t)^{n(1-x)}f(t)dt}{%
\int_{0}^{1}t^{nx}(1-t)^{n(1-x)}dt}
\end{equation*}%
is not stable in the sense of Hyers and Ulam.\newline

\noindent {\large \textbf{(iv) Stancu operators}} (\cite{popa3})\newline

\noindent Let $C[0,1]$ be the linear space of all continuous functions $%
f:[0,1]\rightarrow \mathbb{R}$, endowed with the supremum norm denoted by $%
\Vert .\Vert $, and $a,b$ real numbers, $0\leq a\leq b$. The Stancu operator
\cite{stancu} $S_{n}:C[0,1]\rightarrow \mathit{\Pi _{n}}$ is defined by%
\begin{equation*}
~S_{n}f(x)=\sum\limits_{k=0}^{n}f\biggl{(}\frac{k+a}{n+b}\biggl{)}{\binom{n}{k}}%
x^{k}(1-x)^{n-k},
\end{equation*}%
\ $f\in C[0,1]$. It is HU-stable and the infimum of the Hyers-Ulam constant is $%
K_{S_{n}}={\binom{{2n}}{{2[\frac{n}{2}]}}}/{\binom{{n}}{{[\frac{n}{2}]}}}$,
for each $n\geq 1$.\newline

\noindent {\large \textbf{(v) Kantorovich operators}} (\cite{popa3})\newline

Let $X=\{f;~f:[0,1]\rightarrow \mathbb{R},~%
\mbox{where $f$ is bounded and Riemann
integrable}\}$ be endowed with the supremum norm denoted by $\Vert .\Vert $.
The Kantorovich operator defined by
\begin{equation*}
~K_{n}f(x)=(n+1)\sum\limits_{k=0}^{n}\biggl{(}\int_{\frac{k}{n+1}}^{\frac{k+1%
}{n+1}}f(t)dt\biggl{)}{\binom{n}{k}}x^{k}(1-x)^{n-k},
\end{equation*}%
\newline
$~~f\in X,x\in \lbrack 0,1]$ is stable in Hyers-Ulam sense and the best HUS
constant is $K_{S_{n}}={\binom{{2n}}{{2[\frac{n}{2}]}}}/{\binom{{n}}{{[\frac{%
n}{2}]}}}$.$\vspace{0.5cm}$\newline

\noindent {\large \textbf{3. Main Results}} $\vspace{0.5cm}$\newline

\noindent In this section, we show the Hyers-Ulam stability of some other
operators. Let $D_{R}$ denote the compact disk having radius $R$, i.e., $%
D_{R}=\{z\in \mathbb{C}:~|z|\leq R\}$.\newline

\noindent {\large \textbf{(i) Bernstein-Schurer Operators}} $\vspace{0.5cm}$%
\newline

Let $X_{D_{R}}=\{f:D_{R}\rightarrow \mathbb{C}~\mbox{be analytic in}~D_{R}\}$
be the collection of all analytic functions endowed with the supremum norm
denoted by $\Vert .\Vert $. The complex Bernstein-Schurer operator (\cite{geo})%
\newline
$~~~~~~~~~~~~~S_{n,p}:X_{D_{R}}\rightarrow \it{\Pi _{n+p}}$ is defined by
\begin{equation*}
S_{n,p}(f)(z)=\sum\limits_{k=0}^{n+p}\left(
\begin{array}{c}
n+p \\
k%
\end{array}%
\right) z^{k}(1-z)^{n+p-k}f\biggl{(}\frac{k}{n}\biggl{)},~~z\in \mathbb{C},~f\in X_{D_{R}}.
\end{equation*}%
We have $N(S_{n,p})=\{f\in X_{D_{R}}:~f(\frac{k}{n})=0,~0\leq k\leq n+p\}$,%
~which is a closed subspace of $X_{D_{R}}$, and $R(S_{n,p})=\mathit{\Pi _{n+p}%
}$. The operator
$\widetilde{S}_{n,p}:{X_{D_{R}}/N(S_{n,p})}\rightarrow \mathit{\Pi _{n+p}}$
is bijective, $\widetilde{S}_{n,p}^{-1}:{\mathit{\Pi _{n+p}}}\rightarrow
X_{D_{R}}/N(S_{n,p})$ is bounded since dim$\mathit{\Pi _{n+p}}=2(n+p+1)$. So
according to Theorem 2.2 the operator $S_{n,p}$ is Hyers-Ulam stable.
\newline

\noindent {\large \textbf{Theorem 3.1.}} For $n\geq 1$
\begin{equation*}
K_{S_{n,p}}=\Vert \widetilde{S}_{n,p}^{-1}\Vert ={\binom{2(n+p)}{2[\frac{n+p%
}{2}]}}/{\binom{n+p}{[\frac{n+p}{2}]}}.
\end{equation*}

\noindent {\large \textbf{Proof.}} Let $q\in \mathit{\Pi _{n+p}}$, $\Vert
q\Vert \leq 1$, and its Lorentz representation%
\begin{equation*}
q(z)=\sum\limits_{k=0}^{n+p}c_{k}(q)z^{k}(1-z)^{n+p-k},~~|z|\leq R.\newline
\newline
\end{equation*}%
\newline
Consider the constant function $f_{q}\in X_{D_{R}}$ defined by\newline
\begin{equation*}
f_{q}\biggl{(}\frac{k}{n}\biggl{)}=\frac{c_{k}(q)}{{\binom{n+p}{k}}},~~0\leq k\leq n+p.%
\newline
\newline
\newline
\end{equation*}%
Then $S_{n,p}f_{q}=q$ and $\widetilde{S}_{n,p}^{-1}(q)=f_{q}+N(S_{n,p})$.%
\newline
As usual, the norm of $\widetilde{S}_{n,p}^{-1}:\it{\Pi _{n+p}}\rightarrow {%
X_{D_{R}}}/N(S_{n,p})$ is defined by
\begin{equation*}
\Vert \widetilde{S}_{n,p}^{-1}\Vert =\sup\limits_{\Vert q\Vert \leq 1}\Vert
\widetilde{S}_{n,p}^{-1}(q)\Vert =\sup\limits_{\Vert q\Vert \leq
1}\inf\limits_{h\in N(S_{n,p})}\Vert f_{q}+h\Vert .
\end{equation*}%
Clearly
\begin{equation*}
\inf\limits_{h\in N(S_{n,p})}\Vert f_{q}+h\Vert =\Vert f_{q}\Vert
=\max\limits_{0\leq k\leq n+p}|c_{k}(q)|/{\binom{n+p}{k}}.
\end{equation*}%
Therefore%
\begin{eqnarray*}
\Vert \widetilde{S}_{n,p}^{-1}\Vert &=&\sup\limits_{\Vert q\Vert \leq
1}\max\limits_{0\leq k\leq n+p}|c_{k}(q)|/{\binom{n+p}{k}}\newline
\\
&\leq &\sup\limits_{\Vert q\Vert \leq 1}\max\limits_{0\leq k\leq
n+p}d_{n+p,k}.\Vert q\Vert /{\binom{n+p}{k}}=\max\limits_{0\leq k\leq
n+p}d_{n+p,k}/{\binom{n+p}{k}}.
\end{eqnarray*}%
On the other hand, let $r(z)=T_{n}(2z-1),~|z|\leq R$. Then $\Vert r\Vert =1$
and $|c_{k}(r)|=d_{n+p,k},~0\leq k\leq n+p$, according to Theorem 2.4.
Consequently%
\begin{equation*}
\Vert \widetilde{S}_{n,p}^{-1}\Vert \geq \max\limits_{0\leq k\leq
n+p}|c_{k}(r)|/{\binom{n+p}{k}}=\max\limits_{0\leq k\leq n+p}d_{n+p,k}/{%
\binom{n+p}{k}}\newline
\end{equation*}%
and so%
\begin{equation*}
\Vert \widetilde{S}_{n,p}^{-1}\Vert =\max\limits_{0\leq k\leq n+p}\frac{%
d_{n+p,k}}{{\binom{n+p}{k}}}=\max\limits_{0\leq k\leq n+p}\frac{{\binom{%
2(n+p)}{2k}}}{{\binom{n+p}{k}}}.
\end{equation*}%
Let%
\begin{equation*}
a_{k}=\frac{{\binom{2(n+p)}{2k}}}{{\binom{n+p}{k}}},~0\leq k\leq n+p.
\end{equation*}%
Then%
\begin{equation*}
\frac{a_{k+1}}{a_{k}}=\frac{2n+2p-2k-1}{2k+1},~~0\leq k\leq n+p.
\end{equation*}%
The inequality $\frac{a_{k+1}}{a_{k}}\geq 1$ is satisfied if and only if $%
k\leq \lbrack \frac{n+p-1}{2}]$, therefore%
\begin{equation*}
\max\limits_{0\leq k\leq n+p}a_{k}=a_{[\frac{n+p-1}{2}]+1}=\left\{
\begin{array}{ll}
a_{[\frac{n+p}{2}]}, & \hbox{n+p is even;} \\
a_{[\frac{n+p}{2}]+1}, & \hbox{n+p is odd.}%
\end{array}%
\right.
\end{equation*}%
\newline
Since $a_{[\frac{n+p}{2}]+1}a_{[\frac{n+p}{2}]}$ if $n+p$ is an odd number,
we conclude that%
\begin{equation*}
K_{S_{n,p}}=\Vert \widetilde{S}_{n,p}^{-1}\Vert ={\binom{2(n+p)}{2[\frac{n+p%
}{2}]}}/{\binom{n+p}{[\frac{n+p}{2}]}}.\newline
\end{equation*}

\parindent=8mm This completes the proof of the theorem.\newline

\noindent {\large \textbf{(ii) Kantrovich-Schurer Operators}}\newline
\\
\parindent=8mm Let $X_{D_{R}}=\{f:D_{R}\rightarrow \mathbb{C}~%
\mbox{be analytic
in}~D_{R}\}$ be the collection of all analytic functions endowed with the
supremum norm denoted by $\Vert .\Vert $. The complex Kantrovich-Schurer
operator (\cite{geo})\newline
$~~~~~~~~L_{n,p}:X_{D_{R}}\rightarrow \it{\Pi _{n+p}}$ is defined by%
\begin{equation*}
L_{n,p}(f)(z)=(n+p+1)\sum\limits_{k=0}^{n+p}{\binom{n+p}{k}}%
z^{k}(1-z)^{n+p-k}\int_{\frac{k}{n+1}}^{\frac{k+1}{n+1}}f(t)dt,~~z\in
\mathbb{C},~f\in X_{D_{R}}.
\end{equation*}%
We have%
\begin{equation*}
N(L_{n,p})=\{f\in X_{D_{R}}:~f(t)=0,~~t\in D_{R}\}.
\end{equation*}%
The operators $L_{n,p}$ are Hyers-Ulam stable since their ranges are finite
dimensional spaces.\newline

\noindent {\large \textbf{Theorem 3.2.}} For $n\geq 1$%
\begin{equation*}
K_{L_{n,p}}=\Vert \widetilde{T}_{n,p}^{-1}\Vert =\frac{(n+1){\binom{2(n+p)}{%
2[\frac{n+p}{2}]}}}{(n+p+1){\binom{n+p}{[\frac{n+p}{2}]}}}.
\end{equation*}

\noindent {\large \textbf{Proof.}} Let $q\in \it{\Pi _{n+p}}$, $\Vert q\Vert \leq
1$, and its Lorentz representation%
\begin{equation*}
q(z)=\sum\limits_{k=0}^{n+p}c_{k}(q)z^{k}(1-z)^{n+p-k},~~|z|\leq R.
\end{equation*}

Consider the constant function $f_{q}\in X_{D_{R}}$ defined by%
\begin{equation*}
f_{q}(t)=\frac{(n+1){c_{k}(q)}}{(n+p+1){{\binom{n+p}{k}}}},~~0\leq k\leq
n+p,~~t\in D_{R}.
\end{equation*}%
Then $L_{n,p}f_{q}=q$ and $\widetilde{L}_{n,p}^{-1}(q)=f_{q}+N(L_{n,p})$.%
\newline
As usual, the norm of $\widetilde{L}_{n,p}^{-1}:\it{\Pi _{n+p}}\rightarrow
X_{D_{R}}/N(L_{n,p})$ is defined by%
\begin{equation*}
\Vert \widetilde{L}_{n,p}^{-1}\Vert =\sup\limits_{\Vert q\Vert \leq 1}\Vert
\widetilde{L}_{n,p}^{-1}(q)\Vert =\sup\limits_{\Vert q\Vert \leq
1}\inf\limits_{h\in N(L_{n,p})}\Vert f_{q}+h\Vert .
\end{equation*}%
Clearly%
\begin{equation*}
\inf\limits_{h\in N(L_{n,p})}\Vert f_{q}+h\Vert =\Vert f_{q}\Vert
=\max\limits_{0\leq k\leq n+p}\frac{(n+1)|c_{k}(q)|}{(n+p+1){\binom{n+p}{k}}}%
.
\end{equation*}%
Therefore%
\begin{eqnarray*}
\Vert \widetilde{L}_{n,p}^{-1}\Vert &=&\sup\limits_{\Vert q\Vert \leq
1}\max\limits_{0\leq k\leq n+p}\frac{(n+1)|c_{k}(q)|}{(n+p+1){\binom{n+p}{k}}%
}\newline
\\
&\leq &\sup\limits_{\Vert q\Vert \leq 1}\max\limits_{0\leq k\leq n+p}\frac{%
(n+1)\Vert q\Vert d_{n+p,k}}{(n+p+1){\binom{n+p}{k}}} \\
&=&\max\limits_{0\leq k\leq n+p}\frac{(n+1)d_{n+p,k}}{(n+p+1){\binom{n+p}{k}}%
}.
\end{eqnarray*}%
On the other hand, let $r(z)=T_{n}(2z-1),~|z|\leq R$. Then $\Vert r\Vert =1$
and $|c_{k}(r)|=d_{n+p,k},~0\leq k\leq n+p$, according to Theorem 2.4.
Consequently%
\begin{equation*}
\Vert \widetilde{L}_{n,p}^{-1}\Vert \geq \max\limits_{0\leq k\leq n+p}\frac{%
(n+1)|c_{k}(r)|}{(n+p+1){\binom{n+p}{k}}}=\max\limits_{0\leq k\leq n+p}\frac{%
(n+1)d_{n+p,k}}{(n+p+1){\binom{n+p}{k}}}\newline
\end{equation*}%
and so%
\begin{equation*}
\Vert \widetilde{L}_{n,p}^{-1}\Vert =\max\limits_{0\leq k\leq n+p}\frac{%
(n+1)d_{n+p,k}}{(n+p+1){\binom{n+p}{k}}}=\max\limits_{0\leq k\leq n+p}\frac{%
(n+1){\binom{2(n+p)}{2k}}}{(n+p+1){\binom{n+p}{k}}}.
\end{equation*}%
Let%
\begin{equation*}
a_{k}=\frac{(n+1){\binom{2(n+p)}{2k}}}{(n+p+1){\binom{n+p}{k}}},~0\leq k\leq
n+p.\newline
\newline
\end{equation*}%
Then%
\begin{equation*}
\frac{a_{k+1}}{a_{k}}=\frac{2n+2p-2k-1}{2k+1},~~0\leq k\leq n+p.\newline
\end{equation*}%
The inequality $\frac{a_{k+1}}{a_{k}}\geq 1$ is satisfied if and only if $%
k\leq \lbrack \frac{n+p-1}{2}]$, therefore
\begin{equation*}
\max\limits_{0\leq k\leq n+p}a_{k}=a_{[\frac{n+p-1}{2}]+1}=\left\{
\begin{array}{ll}
a_{[\frac{n+p}{2}]}, & \hbox{n+p is even;} \\
a_{[\frac{n+p}{2}]+1}, & \hbox{n+p is odd.}%
\end{array}%
\right.
\end{equation*}

Since $a_{[\frac{n+p}{2}]+1}=a_{[\frac{n+p}{2}]}$ if $n+p$ is an odd number,
we conclude that
\begin{equation*}
K_{L_{n,p}}=\Vert \widetilde{L}_{n,p}^{-1}\Vert =\frac{(n+1){\binom{2(n+p)}{%
2[\frac{n+p}{2}]}}}{(n+p+1){\binom{n+p}{[\frac{n+p}{2}]}}}.\newline
\end{equation*}

\parindent=8mm This completes the proof of the theorem.\newline

\newpage
\noindent {\large \textbf{(iii) Lorentz Operators}} \newline

The complex Lorentz polynomial \cite{gal3} attached to any analytic function
$f$ in a domain containing the origin is given by
\begin{equation*}
L_{n}(f)(z)=\sum\limits_{k=0}^{n}{\binom{n}{k}}\biggl{(}\frac{z}{n}%
\biggl{)}^{k}f^{(k)}(0),~~n\in \mathbb{N}.
\end{equation*}%
For $R>1$ and denoting $D_{R}=\{z\in \mathbb{C};|z|<R\}$, suppose that $%
f:D_{R}\rightarrow \mathbb{C}$ is analytic in $D_{R}$, i.e., $%
f(z)=\sum\limits_{k=0}^{\infty }c_{k}z^{k}$, for all $z\in D_{R}$.\newline

\noindent {\large \textbf{Theorem 3.3.}} For each $n\geq 1$, the Lorentz
polynomial on compact disk is Hyers-Ulam unstable.\\

\noindent {\large \textbf{Proof.}} Let us denote $e_{j}(z)=z^{j}$, then
from Lorentz operators we can easily obtain that $L_{n}(e_{0})(z)=1$, $L_{n}(e_{1})(z)=e_{1}(z)$ and that for all $%
j,n\in \mathbb{N},~j\geq 2$, we have%
\begin{eqnarray*}
L_{n}(e_{j})(z) &=&{\binom{n}{j}}j!\frac{z^{j}}{n^{j}},~~~~~1\leq R_{1}<R \\
&=&z^{j}\biggl{(}1-\frac{1}{n}\biggl{)}\biggl{(}1-\frac{2}{n}\biggl{)}\cdots \biggl{(}1-\frac{j-1}{n}\biggl{)}.\newline
\end{eqnarray*}%
\newline
Also, since an easy computation shows that\newline
\begin{equation*}
L_{n}(f)(z)=\sum\limits_{j=0}^{\infty }c_{j}L_{n}(e_{j})(z),~~\forall~
|z|\leq R_{1},
\end{equation*}%
and $L_{n}(e_{0})(z)=1$, $L_{n}(e_{1})(z)=e_{1}(z)$. It follows that for
each $j\geq 2,$ $(1-\frac{1}{n})(1-\frac{2}{n})\cdots (1-\frac{j-1}{n})$ is
an eigen value of $L_{n}$. It can be easily seen that $L_{n}$ is injective.
Therefore $1/(1-\frac{1}{n})(1-\frac{2}{n})\cdots (1-\frac{j-1}{n})$ is an
eigen value of $L_{n}^{-1}$. Since%
\begin{equation*}
\lim\limits_{j\rightarrow \infty }\frac{1}{(1-\frac{1}{n})(1-\frac{2}{n}%
)\cdots (1-\frac{j-1}{n})}=\lim\limits_{j\rightarrow \infty }\frac{n^{j}}{%
(n-1)(n-2)\cdots (n-j+1)}=+\infty,
\end{equation*}
we conclude that $L_{n}^{-1}$ is unbounded and so $L_{n}$ is HU-unstable.

\parindent=8mmThis completes the proof of the theorem.


\newpage

\begin{thebibliography}{99}
\bibitem{cam} F. Altomare and M. Campiti, Korovkin-Type Approximation Theory
and its Applications,W. de Gruyter, Berlin, New York, 1994.

\bibitem{aok} T. Aoki, On the stability of linear transformation in Banach
spaces, J. Math. Soc. Japan 2 (1950) 64-66.

\bibitem{geo} G.A. Anastassiou and S.G. Gal, Approximation by complex
Bernstein-Schurer and Kantorovich-Schurer polynomials in compact disks,
Computers and Mathematics with Applications 58 (2009) 734-743.

\bibitem{jung} J. Brzdek and S.M. Jung, A note on stability of an operator
linear equation of the second order, Abstr. Appl. Anal. (2011) 15. Article
ID602713.

\bibitem{ras} J. Brzdek and Th.M. Rassias, Functional Equations in
Mathematical Analysis, Springer, 2011.

\bibitem{gal3} S.G. Gal, Approximation by complex Lorentz polynomials, Math.
Commun., 16 (2011), 67-75.




\bibitem{hat} O. Hatori, K. Kobayasi, T. Miura, H. Takagi and S.E. Takahasi,
On the best constant of Hyers-Ulam stability, J. Nonlinear Convex Anal. 5
(2004) 387-393.

\bibitem{miu} G. Hirasawa and T. Miura, Hyers-Ulam stability of a closed
operator in a Hilbert space, Bull. Korean Math. Soc. 43 (2006) 107-117.

\bibitem{hyer} D.H. Hyers, On the stability of the linear functional
equation, Proc. Natl. Acad. Sci. USA 27 (1941) 222-224.

\bibitem{isac} D.H. Hyers, G. Isac and Th.M. Rassias, Stability of
Functional Equation in Several Variables, Birkh\"{a}user, Basel, 1998.


\bibitem{lorentz} G.G. Lorentz, Bernstein polynomials, 2nd edition,Chelsea
Publ., New York,1986.

\bibitem{lub} D.S. Lubinsky and Z. Ziegler, Coefficients bounds in the
Lorentz representation of a polynomial, Canad. Math. Bull. 33 (1990) 197-206.

\bibitem{lup} A. Lupa\c{s}, Die Folge der Betaoperatoren, Dissertation,
Univ. Stuttgart, 1972.

\bibitem{tak} T. Miura, M. Miyajima and S.E. Takahasi, Hyers-Ulam stability
of linear differential operator with constant coefficients, Math. Nachr. 258
(2003) 90-96.

\bibitem{ansari} M. Mursaleen and K.J. Ansari, Stability results in
intuitionistic fuzzy normed spaces for a cubic functional equation. Appl.
Math. Inform. Sci. 7(5), (2013) 1685-1692.

\bibitem{pol} G. P\'{o}lya and G. Szeg\"{o}, Aufgaben und Lehrs\"{a}tze aus
der Analysis, I, Springer, Berlin, 1925.

\bibitem{rasa2} D. Popa and I. Ra\c{s}a, The Fr\'{e}chet functional equation
with applications to the stability of certain operators, J. Approx. Theory 1
(2012) 138-144.

\bibitem{popa2} D. Popa and I. Ra\c{s}a, On the stability of some classical
operators from approximation theory, Expo. Math. 31(2013) 205-214.

\bibitem{popa3} D. Popa and I. Ra\c{s}a, On the best constant in Hyers-Ulam
stability of some positive linear operators, Jour. Math. Anal. Appl.
412(2014) 103-108.

\bibitem{ras0} Th.M. Rassias, On the stability of the linear mappings in
Banach spaces, Proc. Amer. Math. Soc. 72 (1978) 297-300.

\bibitem{stancu} D.D. Stancu, Asupra unei generaliz\u{a}ri a polinoamelor
lui Bernstein, Stud. Univ. Babe\c{s}-Bolyai 14 (1969) 31-45.

\bibitem{miu1} H. Takagi, T. Miura and S.E. Takahasi, Essential norms and
stability constants of weighted composition operators on C(X), Bull. Korean
Math. Soc. 40 (2003) 583-591.

\bibitem{ulam} S.M. Ulam, A collection of Mathematical problems,
Interscience, New York, 1960.

\bibitem{gajda} Z. Gajda, On stability of additive mappings, Int. J. Math.
Math. Sci. 14 (1991), 431-434.

\bibitem{Brz} J. Brzd\c{e}k, Hyperstability of the Cauchy equation on
resticted domains, Acta Math. Hungar. 141 (2013), 58-67.

\bibitem{bri} N. Brillou\"{e}t-Belluot, J. Brzd\c{e}k, K. Ciepli\'{n}ski,
On some recent developments in Ulam's type stability, Abstr. Appl. Anal.
2012 (2012), Article ID 716936, 41 pp.

\bibitem{rus} I.A. Rus, Remarks on Ulam stability of the operatorial equations,
Fixed Point Theory 10 (2009), 305-320.

\bibitem{urs} C. Urs, Ulam-Hyers stability for coupled fixed points of
cntractive type operators, J. Nonlinear Sci. Appl. 6 (2013), no. 2, 124-136.

\bibitem{sin} W. Sintunavarat, Genaralized Hyers-Ulam stability,well-posedness,
and limit showding of fixed point problems for $\alpha$-$\beta$-contraction
mapping in metric spaces. The Scientific World Journal 2014, Article ID 569174, 7 pp.

\bibitem{rus1} I.A. Rus, Ulam stability of operatorial equations, Functional Equations
in Mathematical Analysis, 287-305, Springer, New York, 2012.

\bibitem{bota} M. Bota, T.p. Petru, G. Petru\c{s}el, Hyers-Ulam stability and
applications in guage spaces. Miskolc Math. Notes 14 (2013), no. 1, 41-47.

\bibitem{bota1} M. Bota, E. Karapinnar, O. Mle\c{s}ni\c{t}e, Ulam-Hyers stability
results for fixed point problems via $\alpha$-$\chi$-contractive mapping in
$(b)$-metric space. Abstr. Appl. Anal. 2013, Art. ID 825293, 6 pp.

\bibitem{pet} A. Petru\c{s}el, G. Petru\c{s}el, C. Urs, Vector-valued metrics,
fixed points and coupled fixed points for nonlinear operators, Fixed Point Theory
Appl. 2013, 2013:218, 21 pp.

\bibitem{brz1} J. Brzd\c{e}k, L. C\u{a}dariu, K. Ciepli\'{n}ski, Fixed Point Theory
and the Ulam stability, J. Function Spaces 2014 (2014), Article ID 829419, 16 pp.

\bibitem{ansari1} S. A. Mohiuddine, M. Mursaleen, Khursheed J. Ansari,
On the Stability of Fuzzy Set-Valued Functional Equations,
The Scientific World Journal, Volume 2014, Article ID 392943,12pages.



\end{thebibliography}
\end{document}